\documentclass[11pt]{amsart}

\usepackage{enumerate}
\usepackage{latexsym}
\usepackage[pdftex]{graphicx}
\usepackage{amssymb}
\usepackage{xcolor}

\newtheorem{theorem}{Theorem}
\newtheorem{lemma}{Lemma}
\newtheorem{proposition}{Proposition}
\newtheorem{remark}{Remark}
\newtheorem{example}{Example}

\newtheorem{definition}{Definition}
\newtheorem{corollary}{Corollary}

\newtheorem{problem}{Problem}

\def\R{{\mathbb R}}

\newcommand{\beq}{\begin{equation}}
\newcommand{\eeq}{\end{equation}}
\newcommand{\beqna}{\begin{eqnarray*}}
\newcommand{\eeqna}{\end{eqnarray*}}
\newcommand{\beqn}{\begin{equation*}}
\newcommand{\eeqn}{\end{equation*}}
\newcommand{\bp}{\begin{proof}}
\newcommand{\ep}{\end{proof}}
\newcommand{\bprop}{\begin{proposition}}
\newcommand{\eprop}{\end{proposition}}
\newcommand{\bt}{\begin{theorem}}
\newcommand{\et}{\end{theorem}}
\newcommand{\bex}{\begin{example}}
\newcommand{\eex}{\end{example}}
\newcommand{\bc}{\begin{corollary}}
\newcommand{\ec}{\end{corollary}}
\newcommand{\bl}{\begin{lemma}}
\newcommand{\el}{\end{lemma}}
\newcommand{\bprob}{\begin{problem}}
\newcommand{\eprob}{\end{problem}}
\newcommand{\br}{\begin{remark}}
\newcommand{\er}{\end{remark}}
\newcommand{\bd}{\begin{definition}}
\newcommand{\ed}{\end{definition}}

\begin{document}

\title
[Directly congruent projections in $\mathbb{R}^5$]
{On   bodies in $\mathbb{R}^5$  with directly congruent projections or sections}

\author[M.A. Alfonseca]{M. Angeles Alfonseca}
\address{Department of Mathematics, North Dakota State University\\
Fargo, ND 58108, USA} \email{maria.alfonseca@ndsu.edu}

\author[M. Cordier]{Michelle Cordier}
\address{Department of Mathematics, Chatham University,
Pittsburgh, PA 15232, USA} \email{m.doyle@chatham.edu}

\author[D. Ryabogin]{Dmitry Ryabogin}
\address{Department of Mathematics, Kent State University,
Kent, OH 44242, USA} \email{ryabogin@math.kent.edu}

\thanks{The  third author is supported in
part by U.S.~National Science Foundation Grant DMS-1600753}

\keywords{Projections and sections of convex bodies}

\begin{abstract}
Let $K$ and $L$  be two convex bodies in ${\mathbb R^5}$ with countably many diameters, such that their projections onto all $4$ dimensional subspaces containing one fixed diameter are directly congruent. We show that if these projections have no rotational symmetries, and the projections of $K,L$ on certain 3 dimensional subspaces have no symmetries, then $K=\pm L$ up to a translation. We also prove the corresponding result for sections of star bodies.
\end{abstract}

\maketitle

\section{Introduction}

In this paper we address the following problems (see \cite[ Problem 3.2, page 125 and Problem 7.3, page 289]{Ga}).
\bprob\label{pr1} 
Suppose that $2\le k\le n-1$ and that $K$ and $L$ are convex bodies in ${\mathbb R}^n$ such that the projection $K|H$ is congruent to $L|H$ for all $H\in {\mathcal G}(n,k)$. Is $K$ a translate of $\pm L$?
\eprob
\bprob\label{pr2}
Suppose that $2\le k\le n-1$ and that $K$ and $L$ are star bodies in ${\mathbb R}^n$ such that the section $K\cap H$ is congruent to $L\cap H$ for all $H\in {\mathcal G}(n,k)$. Is $K$ a translate of $\pm L$?
\eprob
Here we say that $K|H$, the projection of $K$ onto $H$,  is congruent to $L|H$ if there exists   an orthogonal transformation  $\varphi\in O(k,H)$ in $H$ such that $\varphi(K|H)$ is  a translate of  $L|H$; ${\mathcal G}(n,k)$ stands for the Grassmann manifold of all $k$ dimensional subspaces in ${\mathbb R^n}$.

Several partial results are known for Problems \ref{pr1} and \ref{pr2}. For symmetric bodies, the answer is affirmative due to theorems of Aleksandrov (for Problem \ref{pr1}, see \cite{A} and \cite[Theorem 3.3.6, page 115]{Ga}) and Funk (for Problem \ref{pr2}, see \cite[Theorem 7.2.6, page 281]{Ga}). In the class of convex polytopes, the answer to both problems is also affirmative \cite{MyR}. If the projections are translates of each other, or if the bodies are convex and the corresponding sections are translates of each other, again a positive result is obtained (see \cite[Theorems  3.1.3 and 7.1.1]{Ga} and \cite{R1}).  For history and additional partial results, we refer the reader to \cite{ACR}, \cite{My}, \cite{R}, \cite{R2}.

 Hadwiger established that, for $n\geq 4$ and $k=n-1$,  if the orthogonal transformations between the projections are all translations, it is not necessary to know the the projections onto all subspaces; instead, it is enough to have information about the projections on all subspaces containing a fixed line (see \cite{Ha}, and \cite[pages 126--127]{Ga}). 
 
 In this paper, we obtain several Hadwiger-type results for both Problems when $k=4$ in the case of direct congruence; the fixed line will be given by the direction of one of the diameters of the body $K$ (see Section \ref{sec2} for the definitions of direct congruence and diameter).  We follow the ideas from \cite{Go}, \cite{R} and \cite{ACR}, where similar results were obtained in the cases $k=2,3$. The case $k=4$ is harder, due to the fact that four dimensional rotations are more difficult to handle than two or three dimensional ones. Nevertheless, here we obtain the expected conclusion of Problems \ref{pr1} and \ref{pr2} that $K=\pm L$ up to a translation, while in \cite{ACR} the conclusion was that $K=L$ or $K={\mathcal{O}}L$ up to a translation, for a certain orthogonal transformation ${\mathcal{O}}$ of $\mathbb{R}^n$.

\subsection{Results about directly congruent projections}

Let $n\geq 4$ and $S^{n-1}$ be the unit sphere in $\mathbb{R}^n$. Given $w\in S^{n-1}$, let $w^\perp$ be the $(n-1)$ dimensional subspace of $\mathbb{R}^n$ that is orthogonal to $w$. We denote by $d_K(\zeta)$  a diameter of the body $K$ which is parallel to the direction $\zeta \in S^{n-1}$. 

Let $D$ and $B$ be two subsets of $H\in {\mathcal G}(n,k)$, $3 \leq k \leq n-1$. We say that $D$ and $B$ are directly congruent if $\varphi(D)=B+a$ for some vector $a\in H$ and some rotation $\varphi \in SO(k,H)$. We also say that $D$ has an $SO(k)$ symmetry (respectively, $O(k)$ symmetry) if $\varphi(D)=D+a$ for some vector $a\in H$ and some non-identical rotation $\varphi \in SO(k,H)$ (respectively, in $O(k,H)$).

\bigskip

\begin{figure}[ht]
\includegraphics[clip, width=7cm]{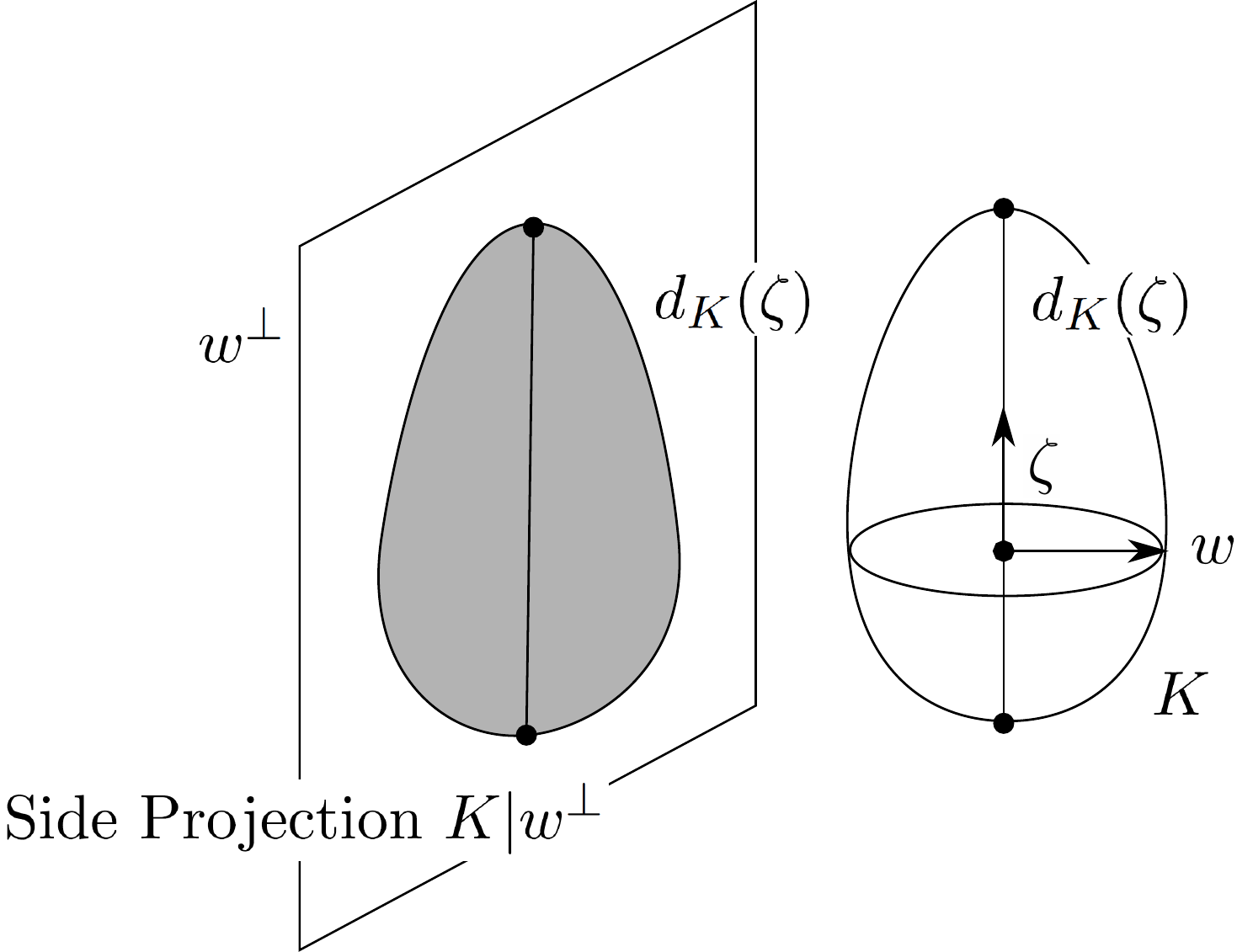}\\
\caption{Diameter $d_K(\zeta)$ and side projection $K|w^\perp$.} 
\label{sgproj}
\end{figure}

We prove the following $5$ dimensional result.
\bt\label{tpr}
Let  $K$ and $L$ be two convex bodies in ${\mathbb R}^5$ having countably many diameters. Assume  that there exists a diameter $d_K(\zeta)$, such that the  side projections $K|H$, $L|H$ onto  all four dimensional subspaces
$H$ containing $\zeta$   are directly congruent, see Figure \ref{sgproj}.
Assume also that these projections  have 
no $SO(4)$ symmetries, and that the three dimensional projections $K|(H \cap \zeta^\perp)$, $L|(H \cap \zeta^\perp)$ have no $O(3)$ symmetries.
Then   $K=L+b$ or $K=-L+b$ for some $b\in {\mathbb R^5}$.

%If, in addition,  the ``ground" projections $K|\zeta^{\perp}$, $L|\zeta^{\perp}$, are directly congruent and do not have rigid motion symmetries, then $K=L+b$ for some $b\in {\mathbb R^4}$.
\et

We state a generalization of Theorem \ref{tpr} to $n$ dimensions as a Corollary.

\bc\label{cpr} Let $K$ and $L$ be two convex bodies in $\mathbb{R}^n$ having countably many diameters.  Assume that there exists a diameter $d_K(\zeta)$, such that the projections $K|H$, $L|H$ onto all four dimensional   subspaces $H$ containing $\zeta$ are directly congruent.  Assume also that these projections have no $SO(4)$ symmetries, and that the three dimensional projections $K|(H\cap \zeta^\perp)$ and $L|(H \cap \zeta^\perp)$ have no $O(3)$ symmetries.  Then $K=L+b$ or $K=-L+b$ for some $b \in \mathbb{R}^n$.
\ec

\subsection{Results about directly congruent sections}  We also obtain results related to Problem \ref{pr2}.

\bt\label{tse}
Let  $K$ and $L$ be two star bodies in ${\mathbb R}^5$ having countably many diameters. Assume  that there exists a diameter $d_K(\zeta)$, containing the origin, such that the  side sections $K \cap H$, $L\cap H$ by  all four dimensional subspaces
$ H$ containing $\zeta$   are directly congruent.
Assume also that these sections  have 
no $SO(4)$ symmetries, and that the three dimensional sections $K\cap (H \cap \zeta^\perp)$, $L\cap (H \cap \zeta^\perp)$ have no $O(3)$ symmetries.
Then   $K=L+b$ or $K=-L+b$ for some $b\in {\mathbb R^5}$ parallel to $\zeta$.
\et

The $n$ dimensional generalization of Theorem \ref{tse} is stated as a Corollary.

\bc\label{cse}
Let  $K$ and $L$ be two star bodies in ${\mathbb R}^n$ having countably many diameters. Assume  that there exists a diameter $d_K(\zeta)$, containing the origin, such that the  sections $K \cap H$, $L\cap H$ by  all four dimensional subspaces $H$ containing $\zeta$   are directly congruent. Assume also that these sections  have 
no $SO(4)$ symmetries, and that the three dimensional sections $K\cap (H\cap \zeta^\perp) $, $L\cap (H \cap \zeta^\perp)$ have no $O(3)$ symmetries.
Then   $K=L+b$ or $K=-L+b$ for some $b\in {\mathbb R^n}$ parallel to $\zeta$.
\ec

The problems still remain open if there is no diameter condition on the bodies. In our proof, this condition allows us to consider only the information on projections or sections on the subspaces containing a diameter, in the spirit of Hadwiger's result. We remark that the class of convex bodies having a finite number of diameters includes all convex polytopes (which are dense in the class of convex bodies with respect to the Hausdorff metric). In fact, the class of convex polytopes whose three and four dimensional projections have no rigid motion symmetries is also dense in the class of convex bodies with respect to the Hausdorff metric (\cite{Pa}, see also \cite[Proposition 2]{ACR}). We further observe that the assumption on the countability of the sets of the diameters of $K$ and $L$ can be weakened (for example, the set of diameters may be taken to be contained in a countable union of  great circles containing $\zeta$). In the proofs below, we only need  that the set of unit vectors parallel to diameters be nowhere dense on the sphere. 

On the other hand, the restriction of a rotation $\Phi \in SO(4)$ to a two dimensional invariant subspace $\Pi$ containing a diameter is an involution ($(\Phi|\Pi)^2=I$), and it seems to be necessary to exclude involutions from our considerations (in \cite{Zh}, a counterexample for Problems 1 and 2  is constructed when the orthogonal transformation is an involution). 

%Involutions in $\mathbb{R}^4$ (which include reflections with respect to a subspace, and rotations by $\pi$) are precisely what we are excluding

The paper is organized as follows. In Section \ref{sec2}, we introduce the needed definitions and  notation. In Section \ref{sec3}, we prove the main auxiliary result of the paper, a functional equation  similar to Proposition 1 in \cite{ACR}.  In Section  \ref{sec4} we prove Theorem \ref{tpr} and Corollary \ref{cpr}, and in Section  \ref{sec5} we prove Theorem \ref{tse} and Corollary \ref{cse}. 

\section{ Notation and auxiliary definitions}
\label{sec2}

\begin{color}{red}\end{color}

We will use the following standard notation. The unit  sphere in ${\mathbb R}^n$, $n\ge 2$, is   $S^{n-1}$. 
Given $w\in S^{n-1}$, the hyperplane orthogonal to $w$ and passing through the origin will be denoted by  $w^{\perp}=\{x\in {\mathbb R^n}:\,x\cdot w=0   \}$, where $x\cdot w=x_1 w_1+\dots+x_n w_n$ is the usual inner product in ${\mathbb R^n}$. The Grassmann manifold of all $k$ dimensional subspaces in ${\mathbb R^n}$ will be denoted by ${\mathcal G}(n,k)$.  The orthogonal group in ${\mathbb R}^n$  is denoted by $O(n)$, and the special orthogonal group in ${\mathbb R}^n$ by $SO(n)$.  If ${\mathcal{U}} \in O(n)$ is an orthogonal matrix, we will write ${\mathcal{U}}^t$ for its transpose. 

We refer to \cite[Chapter 1]{Ga} for the next definitions involving convex and star bodies. A {\em body} in $\R^n$ is a compact set which is equal to the closure of its non-empty interior.  A {\it convex body} is a body $K$ such that for every pair of points in $K$, the segment joining them is contained in $K$.  For $x \in \mathbb{R}^n$, the {\it support function} of a convex body $K$ is defined as 
$
h_K(x)=\max \{x\cdot y:\, \,y\in K   \}
$
 (see page 16 in \cite{Ga}). The {\it width function} $\omega_K(x)$ of $K$ in the direction $x\in S^{n-1}$ is defined as
	$\omega_K(x)=h_K(x)+h_K(-x)$. A segment $[z,y]\subset K$  is called a {\it diameter} of the convex body $K$  if $|z-y|=\max\limits_{\{\theta\in S^{n-1}\}}\omega_K(\theta)$.  We say that a convex body $K\subset {\mathbb R^n}$ has countably many diameters if the width function $\omega_K$ 
reaches its maximum on a countable subset of $S^{n-1}$.

Observe that a convex body $K$ has at most one diameter parallel to a given direction $\zeta \in S^{n-1}$ (for, if  $K$ had two parallel diameters $d_1$, $d_2$, then $K$ would contain a parallelogram with sides $d_1$ and $d_2$, one of whose  diagonals is longer than $d_1$). For this reason, if $K$ has a diameter parallel to $\zeta\in S^{n-1}$, we will denote it by  $d_K(\zeta)$.

A set $S\subset\R^n$ is said to be {\it star-shaped with respect to an interior point $p$}  if the line segment from $p$ to any point in $S$ is contained in $S$.  
For $x\in {\mathbb R}^n\setminus \{0\}$, and $K\subset {\mathbb R}^n$ a nonempty, compact, star-shaped set with respect to the origin, the {\it radial function} of $K$  is defined as 
$
\rho_K(x)=\max \{c:\,cx\in K   \}
$
%Here, the line through $x$ and the origin is assumed to meet $K$
(\cite[page 18]{Ga}). 
%, \begin{color}{red}where in contrast to Gardner's definition, we can consider only non-negative $c$ since we are assuming that $K$ contains the origin)\end{color}. 
We say that a body $K$ is a {\it star body} if  $K$ is star-shaped with respect to the origin and its radial function $\rho_K$ is continuous.

  Given a star body $K$, a segment $[z,y]\subset K$  is called a {\it diameter} of $K$  if $|z-y|=\max\limits_{\{[a,b]\subset K\}}|a-b|$.  If a non-convex star body $K$ has a diameter containing the origin that is parallel to $\zeta\in S^{n-1}$, we will also denote it by $d_K(\zeta)$.

\begin{figure}[ht]
\includegraphics[clip, width=7cm]{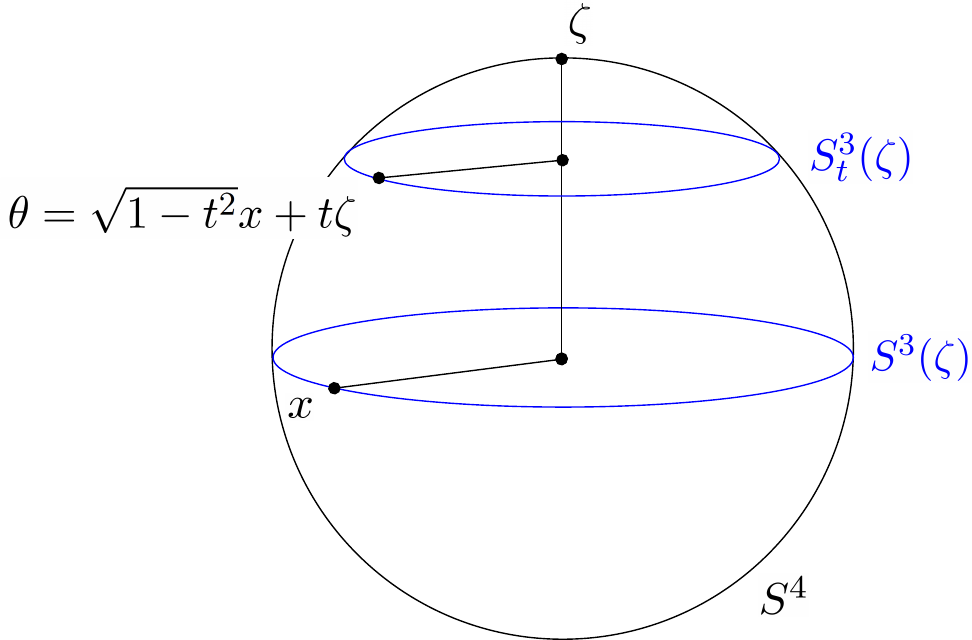}\\
\caption{The great subsphere $S^3(\zeta)$ and the parallel $S^3_t(\zeta)$.} 
\label{parallels}
\end{figure}

 Given $\zeta\in S^{n-1}$,   the great $(n-2)$ dimensional subsphere of $S^{n-1}$ that is perpendicular to $\zeta$ will be denoted by $S^{n-2}(\zeta)=\{\theta\in S^{n-1}:\,\theta\cdot\zeta=0   \}$.   
For $t\in [-1,1]$, the parallel subsphere to $S^{n-2}(\zeta)$ at height $t$ will  be denoted by  $S^{n-2}_t(\zeta)=S^{n-1}\cap \{x\in {\mathbb R^n}:\,x\cdot \zeta=t\}$. Observe that when $t=0$, $S^{n-2}_0(\zeta)=S^{n-2}(\zeta)$. Figure \ref{parallels} shows the case $n=5$.

% Let $E$ be a two or three-dimensional subspace  of ${\mathbb R}^n$.  We will  write $\varphi_{E}\in SO(2, E)$, or $\varphi_{E}\in SO(3, E)$, meaning that  there exists a  choice of an orthonormal basis in ${\mathbb R}^n$ and a rotation $\Phi\in SO(n)$, with a matrix written in this basis, such that  the action of $\Phi$ on $E$ is  the  rotation $\varphi_E$ in $E$, and the  action of $\Phi$ on  $E^{\perp}$ is trivial, {\it i.e.}, $\Phi(y)=y$ for every $y\in E^{\perp}$ (here $E^{\perp}$ stands for the orthogonal complement of $E$). A similar notation will be used for $\varphi_{E}\in O(3, E)$.

For $w  \in S^4$, we will denote by $O(4, S^3(w))$, $SO(4, S^3(w))$, the orthogonal transformations in the $4$ dimensional subspace spanned by the great subsphere $S^3(w)$ of $S^4$. The restriction of a transformation $\varphi \in O(n)$ onto the subspace of smallest dimension containing  $W\subset S^{n-1}$ will be denoted by  $\varphi|_W$. $I$ stands for the identity transformation.

Finally, we define the notion of symmetry for functions, as it will be used throughout the paper.  
% Let $\zeta \in S^4$, and let  $w \in S^3(\zeta)$. For $\alpha \in [0,2]$, we will denote by $\varphi_w^{\alpha\pi}$ the rotation of the sphere $S^2(w)$ by the angle $\alpha\pi$ around $\zeta$, {\it i.e.},  $\varphi_w^{\alpha\pi}(\zeta)=\zeta$. By this we mean that $\varphi_w^{\alpha\pi}$ is the restriction to the  3-dimensional subspace spanned by $S^2(w)$ of a  rotation $\Phi \in SO(4)$  with the following properties:  $\Phi(\zeta)=\zeta$, $\Phi(w)=w$, if  $\{x, y, w, \zeta\}$ is a positively oriented orthonormal basis of  $\mathbb{R}^4$, then for every $v \in (span\{x, y\}\cap S^3)=S^2(w)\cap S^2(\zeta)$, the angle between the vectors $v$ and  $\varphi_w^{\alpha\pi}(v)  \in S^2(w)\cap S^2(\zeta)$  is $\alpha \pi$, and if $\alpha\neq 0,1,2$, $\{v,\varphi_w^{\alpha\pi}(v),w,\zeta\}$ form a positively oriented basis of $\mathbb{R}^4$.

\bd\label{d2}
Let $f$ be a continuous function on $S^{n-1}$ and let $\xi\in S^{n-1}$. We 
 say that the restriction of $f$ onto $S^{k-1}(\xi)$ (or just $f$) has an $SO(k)$ symmetry   if for some non-identical rotation $\varphi_{\xi}\in SO(k, S^{k-1}(\xi))$, we have $f\circ \varphi_{\xi}=f$ on $S^{k-1}(\xi)$. We similarly define the property that $f$ has an $O(k)$ symmetry. 

\ed

\section{A result about a  functional equation on $S^4$}
\label{sec3}

\bprop\label{funeq}
Let $f$ and $g$ be two continuous functions on  $S^4$. Assume that for some $\zeta\in S^4$ and for every $w\in S^3(\zeta)$  there exists a rotation $\varphi_w\in SO(4, S^3(w))$, verifying that 
\begin{equation}
 \label{pos}
f\circ\varphi_w (\theta)=g(\theta), \;\;\; \forall \theta \in S^3(w).
\end{equation}
 Assume, in addition, that  $\varphi_w(\zeta)=\pm \zeta$   $\; \forall w \in S^3(\zeta)$, that the restrictions of $f$ and $g$ to each $S^3(w)$ have no $SO(4)$ symmetries, and that the restrictions of $f$ and $g$ to each $S^3(w)\cap S^3(\zeta)$ have no $O(3)$ symmetries.

Then either $f=g$ on $S^4$ or $f(\theta)=g(-\theta)$  $\forall\theta\in S^4$.
\eprop

\subsection{Auxiliary Lemmata}

We will divide the proof of Proposition \ref{funeq} in several lemmata. The first Lemma describes the structure of the four dimensional rotations $\varphi$ that satisfy the condition $\varphi(\zeta)=\pm\zeta$. Before we state it, let us recall some facts about $SO(4)$ (see, for example, \cite{St} for proofs). 

A rotation  $\varphi \in SO(4)$ such that $\varphi \neq \pm I$  has (at least) two orthogonal invariant two dimensional subspaces, and the restriction of $\varphi$ to each of them is a usual two dimensional rotation. If the two angles corresponding to each of these two subspaces are different, then the subspaces are uniquely determined. In this case, we will denote them by $\Pi$ and $\Pi^\perp$.%, the angle of rotation of $\varphi|\Pi$ by $\alpha$, and the angle of rotation of $\varphi|\Pi^\perp$ by $\beta$. 

\bl\label{rot}
Let $\varphi \in SO(4)$ and let  $\zeta\in \mathbb{R}^4$, $\zeta \neq 0$, such that $\varphi(\zeta)=\pm\zeta$. Assume that $\varphi=\pm I$, then  $\Pi$ and $\Pi^\perp$ are uniquely determined, %(the angles of rotation $\alpha$ and $\beta$ are different),
 without loss of generality  $\zeta \in \Pi^\perp$  and the restriction $\varphi|\Pi^\perp $ is either a trivial rotation or a rotation by  $\pi$.
\el

\bp
Consider first the case in which $\varphi(\zeta)=\zeta$, {\it i.e.} $\zeta$ is an eigenvector of $\varphi$ with eigenvalue $1$. Then $\zeta$ must belong to one of the (possibly many) two dimensional invariant subspaces, which we will call $\Pi^\perp$, and $\varphi|\Pi^\perp $ is a trivial two dimensional rotation. Since $\varphi \neq I$, the orthogonal subspace to $\Pi^\perp$, which we denote by $\Pi$, cannot be an eigenspace associated to $1$. %It must be either the eigenspace associated to $-1$, or the intersection with $\mathbb{R}^4$ of the direct sum of the eigenspaces associated to two complex conjugate eigenvalues. 
Thus, the subspaces $\Pi$ and $\Pi^\perp$ are uniquely determined. The case $\varphi(\zeta)=-\zeta$ can be treated similarly.

\ep

The next Lemma collects some elementary facts about the rotations $\varphi_w$ given by Proposition \ref{funeq}.

\bl\label{controt}
Let $w \in S^3(\zeta)$ and let $\varphi_w \in SO(4,S^3(w))$ be as in the statement of Proposition \ref{funeq}. Then 
\begin{enumerate}[(a)]
\item For each $w\in S^3(\zeta)$, the rotation $\varphi_w$ is uniquely determined.
\item The function mapping $w \rightarrow \varphi_w$ is continuous, and the invariant subspaces $\Pi_w$ and $\Pi_w^\perp$ vary continuously with $w \in S^3(\zeta)$.
\item $\varphi_w$ maps $S^3(\zeta)\cap S^3(w)$ to itself.
\item Let $\mathcal{O}\in SO(5) $ be the orthogonal transformation defined by $\mathcal{O}(\zeta)=\zeta$ and $\mathcal{O}|_{S^3(\zeta)}=-I$. Then  ${\mathcal O}|_{S^3(w)}$ commutes  with $\varphi_w\in SO(4, S^3(w))$.
\end{enumerate}
\el

\bp

\begin{enumerate}[(a)]

\item We show that for each $w\in S^3(\zeta)$, the rotation $\varphi_w \in SO(4,S^3(w))$ given by Proposition \ref{funeq} is unique. On $\mathbb{R}^5$ we consider a positively oriented orthonormal basis $\{u,v,\zeta,z, w\}$, such that $\{u,v,\zeta,z\}$ is a basis for $w^\perp$. We think of $\varphi_w \in SO(4,S^3(w))$ as restriction to $w^\perp$ of a rotation $\Phi_w \in SO(5)$ such that $\Phi_w(w)=w$.

Assume that there are two different rotations $\varphi_1, \varphi_2 \in SO(4,S^3(w))$, verifying 
\[
f\circ\varphi_1 (\theta)=g(\theta), \;\;\; \forall \theta \in S^3(w)
\]
and
\[
f\circ\varphi_2 (\theta)=g(\theta), \;\;\; \forall \theta \in S^3(w).
\]
Then $f\circ\varphi_1  \circ (\varphi_2)^{-1} (\theta)=f(\theta)$ for all $\theta \in S^3(w)$, where $\varphi_1 \circ (\varphi_2)^{-1}$ is not the identity. % since $\varphi_1 \neq \varphi_2$. 
Then $f$ has a rotational symmetry on $S^3(w)$, contradicting the hypotheses of Proposition \ref{funeq}. 

\bigskip

\item We prove the continuity of the map $w \rightarrow \varphi_w$. Let $(w_l)_{l=1}^\infty$ be a sequence of elements of $S^3(\zeta)$ converging to $w \in S^3(\zeta)$ as $l \rightarrow \infty$, and let $\theta$ be any point on $S^3(w)$.  Consider a sequence $(\theta_l)_{l=1}^\infty$ of points $\theta_l \in S^3(w_l)$ converging to $\theta$ as $l \rightarrow \infty$ (to see why such a sequence exists, see \cite[Lemma 3]{ACR}, where an analogous statement is proved). Let $\Phi_{w_l}$ be the rotation in $\mathbb{R}^5$ such that $\Phi_{w_l}(w_l)=w_l$ and  $\Phi_{w_l}|w_l^\perp=\varphi_{w_l}$. By compactness, the sequence $\{\Phi_{w_l}\} \subseteq SO(5)$ has a convergent subsequence. Suppose that $(\Phi_{w_l})_{l=1}^\infty$ has two subsequences that converge to two different rotations in $SO(5)$, $(\Phi_{w_j^1}) \rightarrow \Phi_{w^1}$ and $(\Phi_{w_j^2}) \rightarrow \Phi_{w^2}$ as $j \rightarrow \infty$, where $\Phi_{w^1} \not= \Phi_{w^2}$. Since $w_l$ converges to $w$, and $\Phi_{w_l}(w_l)=w_l$, we have that $\Phi_{w^1}(w)=w,\Phi_{w^2}(w)=w$.
Let $ \varphi_{w^1}$ be the restriction of $\Phi_{w^1}$ to the subspace $w^\perp$, and similarly $\varphi_{w^2}=\Phi_{w^2}|w^\perp$. We know that $f \circ \varphi_{w_j^1} (\theta_j) = g(\theta_j)$ and by passing to the limit as $j \rightarrow \infty$ we obtain that $f \circ \varphi_{w^1} (\theta) = g(\theta)$.  Similarly, we have $f \circ \varphi_{w^2}(\theta) = g(\theta)$.  This implies that $f \circ \varphi_{w^1} (\theta) = f \circ \varphi_{w^2} (\theta)$.  Since the choice of $\theta$ was arbitrary, the last equation holds for all $\theta \in S^3(w)$.  Hence, $f \circ \varphi_{w^1} \circ \varphi_{w^2}^{-1} (\theta) = f (\theta)$ for all $\theta \in S^3(w)$, where $\varphi_{w^1} \circ \varphi_{w^2}^{-1} \not= I$ since $\varphi_{w^1} \not= \varphi_{w^2}$.  Thus, $f$ has a $SO(4)$ symmetry on $S^3(w)$, which is a contradiction.  Therefore, all convergent subsequences of $\{ \Phi_{w_j} \}$ must have the same limit, which we will denote by $\Phi_w$. Since $\Phi_{w_l}(w_l)=w_l$, it follows that  $\Phi_w(w)=w$, and the restriction $\varphi_w=\Phi_w|w^\perp $ is in $SO(4,S^3(w))$. 

If  $\Phi_w|w^\perp \neq \pm I$, let $\Pi_w$ and $\Pi_w^\perp$ be the unique invariant  two dimensional subspaces guaranteed by Lemma \ref{rot}, and assume that $\Phi_w|\Pi_w^\perp=I$ (the case  $\Phi_w|\Pi_w^\perp=-I$ can be treated in a similar way).  Consider a sequence $\left\{\Phi_{w_l} \right\}$ converging to $\Phi_w$. Then, we must have (for a subsequence) that $\Phi_{w_l}|w_l^\perp \neq \pm I$, and letting $\{u_l,v_l,\zeta,z_l,w_l\}$  be the orthonormal basis of $\mathbb{R}^5$  such that $\Pi_{w_l}=span(u_l,v_l)$ and $\Pi_{w_l}^\perp=span(\zeta,z_l)$,  $\Phi_{w_l}|\Pi_{w_l}^\perp=I$. Consider the case where $\Phi_{w_l}|\Pi_{w_l}^\perp= I$ (the other case is analogous). Given that $\Phi_{w_l}(\zeta)=\zeta$ for all $l$, we obtain that $\Phi_w(\zeta)=\zeta$. We also know that $\Phi_{w_l}(z_l)=z_l$ for all $l$, hence there exists $z\in S^3(w)$ such that a subsequence of $\{z_l\}$ converges to $z$, and $\Phi_w(z)=z$. Therefore, $span(\zeta,z)$  must coincide with $\Pi_w^\perp$, and it follows that $\Pi_{w_l}$ converges to $\Pi_w$. 

\bigskip

\item Given $\varphi_w \in SO(4,S^3(w))$ such that $\varphi_w(\zeta)=\pm \zeta$,  by Lemma \ref{rot} there are two invariant two dimensional subspaces, $\Pi_w$ and $\Pi_w^\perp$, where  $\zeta \in \Pi_w^\perp$, and $\varphi_w|\Pi_w^\perp=\pm I$. If $z \in S^3(w)$ is such that $\{\zeta, z\}$ is an orthogonal basis for $\Pi_w^\perp$, we have that $\varphi_w(z)=\pm z$ and hence $span(z)$ is an invariant subspace for $\varphi_w$.

 Since $\zeta$ is orthogonal to $z$ and  $\Pi_w$, we have that $S^3(\zeta)\cap S^3(w)=S^3(w) \cap  \left( \Pi_w \oplus span(z) \right)$. Given that both $\Pi_w$ and  $span(z)$ are invariant subspaces for $\varphi_w$, it follows that $\varphi_w$ maps $S^3(\zeta)\cap S^3(w)$ to itself.

\bigskip

\item  By the definition of ${\mathcal O}$, both $\Pi_w$ and $span(z)$ are eigenspaces for  ${\mathcal O}$, with eigenvalue $-1$, while $span(\zeta)$ is an eigenspace with eigenvalue $1$. Hence,  $\varphi_w$ and ${\mathcal O}$ commute on each of them.  Since $S^3(w)=\Pi_w \oplus span(\zeta) \oplus span(z)$, the result follows.
\end{enumerate}

\ep

\bigskip

\bigskip

The next Lemma is an observation about the geometry of the sphere.

\bl\label{geom}
Let $\zeta$ and $x$ be in $S^k$, $k\ge 3$. Then,
$$
\bigcup\limits_{\{w\in S^{k-1}(\zeta)\cap S^{k-1}(x)  \}}S^{k-1}(w)=S^k.
$$ 
\el
\bp
Let $y$ be any point on $S^k$. Then $S^{k-1}(\zeta)\cap S^{k-1}(x)\cap S^{k-1}(y)$ is nonempty, since $k \geq 3$. Taking any  $w\in S^{k-1}(\zeta)\cap S^{k-1}(x)\cap S^{k-1}(y)\subset S^{k-1}(\zeta)\cap S^{k-1}(x)$, it follows that $y\in S^{k-1}(w)$.
\ep

\bl\label{symmetry}(cf. Lemma 1, \cite{R}).
Let $\zeta\in S^k$, $k\ge 4$. If for every  $w\in S^{k-1}(\zeta)$ we have either $f(\theta)=g(\theta)$ for all $\theta\in S^{k-1}(w)$ or $f(-\theta)=g(\theta)$ for all $\theta\in S^{k-1}(w)$, then 
either $f=g$ on $S^k$ or  $f(-\theta)=g(\theta)$ for all $\theta\in S^k$.
\el
\bp
Assume at first that  there exists  an $x\in S^k$ such that for all $w\in S^{k-1}(\zeta)\cap S^{k-1}(x)$ we have $f(\theta)=g(\theta)$ for all $\theta\in S^{k-1}(w)$.
Then, using  the previous lemma,  we obtain $f=g$ on $S^k$.

Assume now that  there exists  an $x\in S^k$ such that for all $w\in S^{k-1}(\zeta)\cap S^{k-1}(x)$ we have $f(-\theta)=g(\theta)$ for all $\theta\in S^{k-1}(w)$. Then, using the previous lemma, we obtain $f(-\theta)=g(\theta)$ for all $\theta \in S^k$.

Finally, assume that for every $x\in S^k$ there exist two directions $w_1$ and $w_2$ in $S^{k-1}(\zeta) \cap S^{k-1}(x)$ such that $f(\theta)=g(\theta)$ for all $\theta \in S^{k-1}(w_1)$ and $f(-\theta)=g(\theta)$ for all $\theta \in S^{k-1}(w_2)$. Then $f(-x)=f(x)=g(x)$, and since $x$ was chosen arbitrarily, we obtain $f=g$ on $S^k$.
\ep

\bigskip

Let $\mathcal{O}\in SO(5) $ be the orthogonal transformation defined in Lemma \ref{controt} (d)  by $\mathcal{O}(\zeta)=\zeta$ and $\mathcal{O}|_{S^3(\zeta)}=-I$. %Observe  that ${\mathcal O}|_{S^3(w)}$ commutes with every rotation $\varphi_w\in SO(4, S^3(w))$, such that  $\varphi_w(\zeta)= \pm \zeta$, where $w\in S^3(\zeta)$.
A function $f$ defined on $S^4$ can be decomposed in the form
\begin{equation}\label{nax2}
f(\theta)=\frac{f(\theta)+f(\mathcal O\theta)}{2}+\frac{f(\theta)-f(\mathcal O\theta)}{2}=f_{\mathcal O, e}(\theta)+f_{\mathcal O, o}(\theta),\quad\theta\in S^4,
\end{equation}
and we will call $f_{\mathcal O, e}$, $f_{\mathcal O, o}$, the even and odd parts of $f$ with respect to ${\mathcal O}$.
Since ${\mathcal O}^2=I$,  we have
$$
f_{\mathcal O, e}(\theta)=f_{\mathcal O, e}(\mathcal O\theta),\qquad f_{\mathcal O, o}(\theta)=-f_{\mathcal O, o}(\mathcal O\theta).
$$

Given $y\in S^4$, we have that $y \in S^3_t(\zeta)$ for some $t \in [-1,1]$, {\it i.e.} we can write 
\begin{equation}\label{fru1}
y=\sqrt{1-t^2}x+t\zeta,
\end{equation} 
for some $t\in [-1,1]$ and $x\in S^3(\zeta)$ (see Figure \ref{parallels}). For any function $f$   on $S^4$, we define the  function $F_t$  on $S^3(\zeta)$,
\begin{equation}\label{ara}
F_t(x)=F_{t,\zeta}(x)=f(\sqrt{1-t^2}x+t\zeta),\qquad x\in S^3(\zeta),
\end{equation}
which is the restriction of $f$ to $S^3_t(\zeta)$.
Observe that the even part of $F_t$, $(F_t)_e$ equals
$$
{(F_t)}_e(x)=\frac{f(\sqrt{1-t^2}x+t\zeta)+f(-\sqrt{1-t^2}x+t\zeta)}{2}=\frac{f(y)+f(\mathcal Oy)}{2},
$$
where $y$ is as in (\ref{fru1}), {\it  i.e.},
\begin{equation}\label{uru1}
{(F_t)}_e(x)=f_{\mathcal O, e}(y),\qquad {(F_t)}_o(x)=f_{\mathcal O, o}(y).
\end{equation}
Note that  $(F_t)_e(x)=(F_t)_e(-x)$ for every $x\in S^3(\zeta)$. We similarly define $G_t$ from the function $g$.

Every two dimensional great subpshere of $S^3(\zeta)$ is of the form $E_w:=S^3(w)\cap S^3(\zeta)$ for some $w\in S^3(\zeta)$. For $\varphi_w   \in SO(4,S^3(w))$ as in Proposition \ref{funeq}, we have that $\varphi_w(E_w)=E_w$ by Lemma \ref{controt} (c). Denote by  $\phi_{E_w}=\varphi_w|E_w$ the restriction of $\varphi_w$ to $E_w$, and define $G_t$ from $g$ similarly to $F_t$ in (\ref{ara}). Observe  that if $y$ is as in \eqref{fru1} with $x\in E_w$, then $\varphi_w(y)=\sqrt{1-t^2}\phi_{E_w}(x)\pm t\zeta$, and it follows from \eqref{pos}  that, for every $t\in [-1,1]$, 
 \begin{equation}\label{ar111}
F_t\circ\phi_{E_w}(x)=G_t(x)\qquad\forall x\in E_w.
\end{equation}

 \bigskip

\bl\label{radon}
 Assume that  $f,g$ satisfy equation \eqref{pos} for all $w \in S^3(\zeta)$. Then $f_{\mathcal{O},e}(y)=g_{\mathcal{O},e}(y)$ for every $y \in S^4$.
 \el
\bp
For $w \in S^3(\zeta)$, we consider the spherical Radon transform 
\[
   Rf(w)=\int_{E_w}f(x) dx
\]
(see \cite[pg. 429]{Ga}). Since Lebesgue measure is invariant under orthogonal transformations on $E_w$, by \eqref{ar111} we have 
\[
    \int_{E_w} F_t(x)dx=\int_{E_w} F_t(\varphi_w(x))dx=\int_{E_w} G_t(x)dx, 
\]
for each $t\in [-1,1]$.
Hence, $RF_t(w)=RG_t(w)$ for every $w \in S^3(\zeta)$ and $t\in [-1,1]$, and it follows from Proposition 3.4.12  \cite[pg. 108]{Gr} %Theorem C.2.4 from \cite[pg. 430]{Ga} 
that the even parts of $F_t$ and $G_t$ coincide. By equations \eqref{ara} and \eqref{uru1} and Lemma \ref{controt} (d), this means that $f_{\mathcal{O},e}(y)=g_{\mathcal{O},e}(y)$ for every $y \in S^4$.\ep

%\bigskip

{\bf  Note:} Because of Lemma \ref{radon}, from now on we will assume that $f,g$ are odd with respect to $\mathcal{O}$. 

\bigskip

%\begin{color}{red}Observe that it is possible to have a rotation by angle $0$ or $\pi$ on $\Pi$, even under the assumption that $\varphi_w \neq \pm I$. This is the case if $\varphi_w|\Pi=I$ and $\varphi_w|\Pi^\perp=-I$ or viceversa. \end{color}

By Lemmas \ref{rot} and  \ref{controt}, given $w\in S^3(\zeta)$ and  $\varphi_w \neq \pm I$ verifying the hypotheses of Proposition \ref{funeq}, we have that $w^\perp=\Pi_w \oplus \Pi_w^\perp$, where $\zeta \in \Pi_w^\perp$, $\varphi_w |\Pi_w^\perp=\pm I$, and   
$\varphi_w|\Pi_w$ is a 2 dimensional rotation.  %On $\mathbb{R}^5$ we consider a positively oriented orthonormal basis $\{u,v,\zeta,z, w\}$, so that $\{u,v\}$ is a basis of $\Pi_w$, and $\{\zeta, z \}$ is a basis of $\Pi_w^\perp$. We think of $\varphi_w \in SO(4,S^3(w))$ as restriction to $w^\perp$ of a rotation $\Phi \in SO(5)$ with the following properties:  $\Phi(w)=w$, $\Phi|\Pi_w^\perp=\pm I$, and $\Phi|\Pi_w$ is a 2 dimensional rotation. 
Given $t \in \Pi_w \cap S^4$, if the angle between the vectors $t$ and $\varphi_w(t) \in \Pi_w \cap S^4$ is $\alpha \pi$, for $\alpha \in (0,2)$, $\alpha  \not=1$, and $\{t, \varphi_w(t), \zeta, z, w \}$ forms a positively oriented basis of $\mathbb{R}^5$, then  we will denote $\varphi_w|\Pi_w$ by $\varphi_w^{\alpha \pi}$ when we want to specify the angle of rotation.

\bigskip

 We define the sets
\[
   \Xi_+=\left\{w\in S^3(\zeta):\,\ f (\theta)=g(\theta) \,\,\forall \theta\in S^3(w) \right\},
\]
\[
   \Xi_-=\left\{w\in S^3(\zeta):\,\ f (\theta)=g(-\theta) \,\,\forall \theta\in S^3(w) \right\},
\]
\[
   \Xi_0=\left\{w\in S^3(\zeta):\,\  %w^\perp=\Pi_w \oplus \Pi_w^\perp, 
   f|\Pi_w=g|\Pi_w, \textrm{ and } f(\theta)=g(-\theta) \,\,\forall \theta\in \Pi_w^\perp \right\},
\]
\[
   \Xi_1=\left\{w\in S^3(\zeta):\,\  %w^\perp=\Pi_w \oplus \Pi_w^\perp, 
   f|\Pi_w^\perp=g|\Pi_w^\perp, \textrm{ and } f(\theta)=g(-\theta) \,\,\forall \theta\in \Pi_w \right\},
\]
and, for $\alpha\in (0,1)\cup (1,2)$,   
\begin{equation}\label{vasya2}
\Xi_{\alpha}=\left\{w\in S^3(\zeta):\,%w^\perp=\Pi_w \oplus \Pi_w^\perp, 
\varphi_w|\Pi_w^\perp=\pm I, \varphi_w|\Pi_w=\varphi_w^{\alpha \pi} \right.
\end{equation}
\[ \left. \textrm{ and } f\circ\varphi_w (\theta)=g(\theta), \;\;\; \forall \theta \in S^3(w) \right\}.
\]

With this notation, the hypothesis of Proposition \ref{funeq} is that 
\begin{equation}
  \label{baire0}
S^3(\zeta)=\Xi_+ \cup \Xi_- \cup \bigcup_{\alpha \in [0,2)} \Xi_\alpha,
\end{equation}
and we want to show that, under the condition on the lack of symmetries, we have either $S^3(\zeta)=\Xi_+$ or  $S^3(\zeta)=\Xi_-$. By Lemma \ref{symmetry}, this will imply that either $f=g$ on $S^4$ or $f(\theta)=g(-\theta)$ for all $\theta\in S^4$.

We observe that the sets on the right hand side of \eqref{baire} are  disjoint. Indeed, assume that  there exists $z \in \Xi_\alpha \cap \Xi_\beta$,  where  $\alpha$ and $\beta$ are either $+,-$ or a rational number in $[0,2)$, and $\alpha \neq \beta$. Then, there are two rotations $\varphi_{z,\alpha}, \varphi_{z,\beta} \in SO(4,S^3(z))$ such that $f\circ\varphi_{z,\alpha} (\theta)=g(\theta)$ and  $f\circ\varphi_{z,\beta} (\theta)=g(\theta)$ for all $\theta \in S^3(z)$. But then $f\circ\varphi_{z,\beta} (\theta) \circ (\varphi_{z,\alpha})^{-1} (\theta)=f(\theta)$ for all $\theta \in S^3(z)$, where $\varphi_{z,\beta} (\theta) \circ (\varphi_{z,\alpha})^{-1}$ is not the identity since $\alpha \neq \beta$. Thus, if $\Xi_\alpha \cap \Xi_\beta \neq \emptyset$,  $f$ has a symmetry on $S^3(z)$, contradicting the hypothesis of Proposition \ref{funeq}.

\bl\label{closed} The sets $\Xi_+,\Xi_-,\Xi_\alpha$ are closed.
\el
\bp
Since the empty set is closed, we can assume that the sets $\Xi_+,\Xi_-$ and $\Xi_\alpha$ are not empty. First we prove that $\Xi_+$ is closed.  Let $(w_l)_{l=1}^\infty$ be a sequence of elements of $\Xi_+$ converging to $w \in S^3(\zeta)$ as $l \rightarrow \infty$, and let $\theta$ be any point on $S^3(w)$.  Consider a sequence $(\theta_l)_{l=1}^\infty$ of points $\theta_l \in S^3(w_l)$ converging to $\theta$ as $l \rightarrow \infty$.  %(To see why such a sequence exists, see \cite[Lemma 3]{ACR}). 
By definition of $\Xi_+$ we have the following,
\begin{equation*}\label{Xi_plus}
f(\theta_l)=g(\theta_l), \,\,\, \theta_l \in S^3(w_l), \,\,\, l\in \mathbb{N}.
\end{equation*}
Since $f$ and $g$ are continuous, we may pass to the limit and obtain $f(\theta)=g(\theta)$.  Thus $w\in \Xi_+$ since the choice of $\theta \in S^3(w)$ was arbitrary, and hence $\Xi_+$ is closed. A similar proof shows that $\Xi_-$ is closed, replacing $g(\theta_l)$ with $g(-\theta_l)$ and $g(\theta)$ with $g(-\theta)$.

Now we prove that $\Xi_\alpha$ is closed, where $\alpha \in [0,2)$.  As above, let $(w_l)_{l=1}^\infty$ be a sequence of elements of $\Xi_\alpha$ converging to $w \in S^3(\zeta)$ as $l \rightarrow \infty$. %as above,  let $\theta$ be any point on $S^3(w)$ and consider a sequence $(\theta_l)_{l=1}^\infty$ of points $\theta_l \in S^3(w_l)$ converging to $\theta$ as $l \rightarrow \infty$.  
 Let $\varphi_{w_l}$ be the rotation associated to $w_l^\perp$ by Proposition \ref{funeq}. By Lemma \ref{controt} (b), the limit of $\varphi_{w_l}$ is the rotation $\varphi_w$, and if 
$\Pi_{w_l}$ and $\Pi_{w_l}^\perp$ are the two dimensional invariant subspaces of $\varphi_{w_l}$, and  $\Pi_w$ and $\Pi_w^\perp$ are those of $\varphi_w$, we have $(\Pi_{w_l}) \rightarrow \Pi_w$ and $(\Pi_{w_l}^\perp) \rightarrow \Pi_w^\perp$ as $l \rightarrow \infty$.  
Hence, there is either a subsequence $\varphi_{w_j^1}|\Pi_{w_j^1}^\perp=I$ (which would imply $\varphi_w|\Pi_w^\perp=I$), or a subsequence $\varphi_{w_j^1}|\Pi_{w_j^1}^\perp=-I$ (which implies that $\varphi_w|\Pi_w^\perp=-I$). Furthermore, since all $\varphi_{w_l}|\Pi_{w_l}$ are rotations by the angle $\alpha \pi$, for the limit we obtain that $\varphi_w|\Pi_w$ is also a rotation by the same angle, and since
\begin{equation}\label{w_l}
f \circ \varphi_{w_l}^{\alpha \pi}(\theta_l)=g(\theta_l) \;\;\; \theta_l \in \Pi_{w_l}, \;\;\; l \in \mathbb{N}, 
\end{equation}
we conclude by continuity that  $f \circ \varphi_w(\theta)=g(\theta)$. This shows that $\Xi_{\alpha}$ is closed.

\ep

%\bigskip

The next Lemma shows that rotations by an irrational multiple of $\pi$ do not occur because of the lack of symmetries of $f$ and $g$. 

\bl\label{irra}
Under the hypotheses of Proposition \ref{funeq}, we have $\Xi_{\alpha}=\emptyset$ for $\alpha\in \left({\mathbb R}\setminus{\mathbb Q}\right) \cap [0,2)$.
\el
\bp

%The proof of this Lemma is similar to that of Lemma 4 in \cite{ACR}. 
Let $\alpha  \in \left({\mathbb R}\setminus{\mathbb Q}\right) \cap [0,2)$, and take $w \in \Xi_\alpha$. 
 Following the ideas of Schneider \cite{Sch1}, we claim at first that $f^2=g^2$ on $S^3(w)$.
Indeed, since $f$ and $g$ are odd with respect to $\mathcal{O}$, $f^2$ and $g^2$ are even with respect to $\mathcal{O}$, and  \eqref{pos} holds with $f^2$, $g^2$ instead of $f$, $g$. Thus, by Lemma \ref{radon}, we obtain        that  $f^2=g^2$ on $S^3(w)$. 

Squaring (\ref{pos}), we have %(with $\varphi_w=\varphi_w^{\alpha \pi}$), 
$$
f^2\circ\varphi_w(\theta)=g^2(\theta)=f^2(\theta)\qquad \forall\theta\in S^3(w).
$$
Iterating, for any $k\in{\mathbb Z}$,
\begin{equation}\label{iterate}
f^2\circ\varphi^k_w(\theta)=f^2\circ\varphi^{k-1}_w(\theta)=\dots=f^2(\theta)\qquad \forall \theta\in S^3(w).
\end{equation}

Let $\{\zeta,z\}$ be an orthonormal basis for $\Pi^\perp_w$, and consider the three dimensional subspace  generated by $\Pi_w$ and $z$, and  its unit sphere $S^3(w) \cap S^3(\zeta)$. If $\varphi_w|\Pi^\perp_w=I$, then $\varphi_w|S^3(w) \cap S^3(\zeta)$ is a rotation of angle $\alpha \pi$ around the vector $z$. For each $\theta \in S^3(w) \cap S^3(\zeta)$, equation \eqref{iterate} holds for any $k\in \mathbb{Z}$, and since the orbit of $(\varphi^k_w(\theta))_{k\in{\mathbb Z}}$ is dense, we conclude that $f^2$ and $g^2$ are constant on each parallel of $S^3(w) \cap S^3(\zeta)$ perpendicular to $z$.  By continuity, $f$ and $g$ must also be constant on each parallel, and thus $f\circ\varphi_w(\theta)=f(\theta)$ for every $\theta \in S^3(w) \cap S^3(\zeta)$. But then $f$ has an $SO(3)$ symmetry on 
$S^3(w) \cap S^3(\zeta)$, contradicting the hypothesis of Proposition \ref{funeq}.

On the other hand, if $\varphi_w(\zeta)=-\zeta$, then $\varphi_w^2|S^3(w) \cap S^3(\zeta)$ is a rotation around $z$ by the angle $2 \alpha \pi$, and similarly to the previous case, $f$ must be constant on every parallel of $S^3(w) \cap S^3(\zeta)$, and thus has a rotational symmetry on $S^3(w) \cap S^3(\zeta)$. This is a contradiction, and thus we have proven that  $\Xi_\alpha=\emptyset$ for irrational $\alpha$.
\ep

%\begin{color}{red}It could be possible for $S^3(\zeta)$ to be the union of two of the $\Xi$ sets with different angle without contradicting the lack of symmetries if the angles are multiples of each other, for example $\pi/4$ and $\pi/2$. this will not be the case for  incomparable angles, but we have to be careful about this case.
%We also need to make sure that all angles are the same. If they are allowed to vary from subspace to subspace, they could eventually become zero, and then we would not have closed sets. So it is crucial for our purposes that $S^3(\zeta)$ is just one of the sets. 
%\end{color} 

{\bf Proof of Proposition \ref{funeq}.} By \eqref{baire0} and Lemma \ref{irra}, we have written $S^3(\zeta)$ as a countable union of disjoint closed sets,
\begin{equation}
  \label{baire}
   S^3(\zeta)=\Xi_+ \cup \Xi_- \cup \bigcup_{\alpha \in [0,2)\cap \mathbb{Q}} \Xi_\alpha.
\end{equation}

By a well-known result of Sierpi\'nski's \cite{Si}, $S^3(\zeta)$ must equal just one of the sets. We will now assume that $S^3(\zeta)=\Xi_\alpha$ for some $\alpha \in [0,2)\cap \mathbb{Q}$, and derive a contradiction. This will leave us only with the possibilities  $S^3(\zeta)=\Xi_+ $ or $S^3(\zeta)=\Xi_-$, and Proposition \ref{funeq} will be proven.

\bigskip

Let $S^3(\zeta)=\Xi_\alpha$. Choose $w \in S^3(\zeta)$ and $\xi \in S^3(w) \cap S^3(\zeta)$. By Lemma \ref{rot},  the subspace $\xi^\perp$ is equal to $\Pi_\xi \oplus \Pi_\xi^\perp$, with $\zeta \in \Pi_\xi^\perp$. Let $z_\xi \in \Pi_\xi^\perp$ be such that $\{\zeta,z_\xi\}$ is an orthonormal basis of $\Pi_\xi^\perp$. 

We will write the set $S^3(w) \cap S^3(\zeta)$ as the union of two closed sets, $\Theta_{good}$ and $\Theta_{bad}$, where 
\[
\Theta_{good}=\{ \xi \in S^3(w) \cap S^3(\zeta):  \, \Pi_\xi = \xi^\perp \cap w^\perp \cap \zeta^\perp   \},
\]
\[
\Theta_{bad}=\{ \xi \in S^3(w) \cap S^3(\zeta):  \, \textrm{dim}(\Pi_{\xi}\cap w^{\perp} \cap \zeta^{\perp})=\textrm{dim}(\Pi_{\xi}^{\perp}\cap w^{\perp} \cap \zeta^{\perp})=1\}.
\]
Observe that $\Theta_{good}$ and $\Theta_{bad}$ are  closed. Indeed, for any $\xi \in \Theta_{good}$,    $\Pi_\xi^\perp$ is spanned by $\{\zeta,w\}$, and the result follows from Lemma \ref{controt} (b). 
If $\xi \in \Theta_{bad}$, then  $\Pi_\xi^\perp$ is spanned by $\{\zeta,z_\xi\}$ and $\Pi_\xi$ is spanned by $\{w,v_\xi\}$ for some $v_\xi$. Let $\{\xi_j\}$ be a sequence from $\Theta_{bad}$ converging to $\xi \in S^3(w)\cap S^3(\zeta)$. Applying Lemma \ref{controt} (b) once again, we obtain that $\Pi_{\xi_j}$ converges to $\Pi_\xi$  and $\Pi_{\xi_j}^\perp$ converges to $\Pi_\xi^\perp$ as $j \rightarrow \infty$. In particular, $w \in \Pi_\xi$, while $\zeta \in \Pi_\xi^\perp$, and thus $\xi \in \Theta_{bad}$.

%. Passing to a subsequence if necessary, let $z=\lim z_{\xi_j}$ and $v= \lim v_{\xi_j}$. Since $\varphi_{\xi_j}(z_{\xi_j})=\pm z_{\xi_j}$, it follows by continuity (Lemma \ref{controt} (b)) that $\varphi_\xi (z)=\pm z$, and $\{\zeta, z\}$ is an invariant two-dimensional subspace for $\varphi_\xi$ with eigenvalue $\pm 1$. On the other hand, $\varphi_{\xi_j}|\Pi_{\xi_j}$ is a rotation by the angle $\alpha \pi$, and therefore so is $\varphi_{\xi}|\Pi_\xi$. Since $\varphi_{\xi_j}(w)$ is a vector in $\Pi_{\xi_j}$ forming an angle $\alpha \pi$ with $w$, the same holds for $\varphi_\xi(w)$ and $w$. This means that $w \in \Pi_\xi$, while $\zeta \in \Pi_\xi^\perp$, and thus $\xi \in \Theta_{bad}$.

The sets  $\Theta_{good}$ and $\Theta_{bad}$  are {\it disjoint} and their union equals $S^3(w) \cap S^3(\zeta)$. It follows that either $S^3(w) \cap S^3(\zeta)=\Theta_{good}$ or $S^3(w) \cap S^3(\zeta)=\Theta_{bad}$.

Assume that  $S^3(w) \cap S^3(\zeta)=\Theta_{bad}$. We claim that, in this case, the map 
 \[
\xi\to  \ell(\xi)=\Pi_{\xi}^\perp \cap w^{\perp} \cap \zeta^{\perp}
\]
defines a non-vanishing continuous tangent line field on the two dimensional sphere  $S^3(w) \cap S^3(\zeta)$. 
%In order to see this, let $\xi \in S^3(w) \cap S^3(\zeta)$. For $\xi^\perp=\Pi_\xi \oplus \Pi_\xi^\perp$, let $\{u_\xi, v_\xi \}$ be an orthonormal basis of $\Pi_\xi$, such that $\Pi_{\xi}\cap \zeta^{\perp}\cap w^{\perp}=span(u_\xi)$, and recall that  $\{\zeta,y_\xi\}$ is an orthonormal basis of $\Pi_\xi^\perp$. We will show that, as $\xi$ varies continuously on the two dimensional sphere $S^3(w)\cap S^3(\zeta)$, the line spanned by $y_\xi$ also moves continuously. 
If this map were not continuous, then there would exist two subsequences $\{\xi^1_j\}$ and $\{\xi^2_j\}$, both with limit $\xi_0$, such that
\[
 \lim_{j \rightarrow \infty} \ell(\xi^1_j) \neq \lim_{j \rightarrow \infty} \ell(\xi^2_j).
\]
Denote by $z_0^1$ a unit vector in the direction of the line $ \lim_{j \rightarrow \infty} \ell(\xi^1_j)$ and by $z_0^2$ a unit vector in the direction of the line $ \lim_{j \rightarrow \infty} \ell(\xi^2_j)$. We have $z_0^1 \neq \pm z_0^2$. Let $\varphi_{\xi_0}^i=\displaystyle \lim_{j \rightarrow \infty} \varphi_{\xi_j^i}$, for $i=1,2$, and let $\Pi_0^i \oplus (\Pi_0^i)^\perp$ the corresponding decompositions of $\xi_0^\perp$. Since 
%in the proof of Lemma \ref{closed},  since 
all the rotations $\varphi_{\xi_j^i}|\Pi_{\xi_j^i}$
are by the angle $\alpha \pi$, the limiting  rotations $\varphi_{\xi_0}^1|\Pi_0^1$ and $\varphi_{\xi_0}^2|\Pi_0^2$ are by the angle $\alpha \pi$ as well (this can be shown by a reasoning similar to the one in the proof of Lemma \ref{closed}). But given that $z_0^1 \neq \pm z_0^2$, the corresponding limiting subspaces $(\Pi_0^1)^\perp$ and $(\Pi_0^2)^\perp$ 
must be different. This means that  $\varphi_{\xi_0}^1$ and $\varphi_{\xi_0}^2$ are two different rotations on $\xi_0^\perp$. From equation \eqref{pos}, it follows that  $f \circ \varphi_{\xi_0}^1 \circ (\varphi_{\xi_0}^2)^{-1} =f$ on $S^3(\xi_0)$, where $\varphi_{\xi_0}^1 \circ (\varphi_{\xi_0}^2)^{-1} \neq I$. Therefore,
we conclude that $f$ has a rotational symmetry on $\xi_0^\perp$. Thus, the line field spanned by $\ell(\xi)$ must be continuous on the two dimensional sphere $S^3(w) \cap S^3(\zeta)$. This is impossible by a well known result of Hopf (see \cite{Sa}).

\bigskip

We now consider the case in which $S^3(w) \cap S^3(\zeta)=\Theta_{good}$, {\it i.e.} the two dimensional space $\Pi_\xi$ is equal to $\xi^\perp \cap w^\perp \cap \zeta^\perp $, and the restriction of $\varphi _\xi$ to this subspace, which we will denote by $\psi_\xi$, is a rotation by the angle $\alpha \pi$. We have that the restrictions of $f$ and $g$ to $S^3(w) \cap S^3(\zeta)$ satisfy
\begin{equation}
  \label{restr}
  f \circ \psi_\xi (\theta)=g(\theta),  \; \forall \theta \in \Pi_\xi,  
\end{equation}
for every $\xi \in S^3(w) \cap S^3(\zeta)$.  But every one dimensional great circle on the  two dimensional sphere $S^3(w) \cap S^3(\zeta)$ is of the form $S^3(\xi) \cap S^3(w) \cap S^3(\zeta)$ for some $\xi \in S^3(w) \cap S^3(\zeta)$. We are thus under the hypothesis of the continuous Rubik's cube \cite{R}, and therefore we can conclude that either $f=g$ on $S^3 (w) \cap S^3(\zeta)$, or  $f(\theta)=g(-\theta)$ for every $\theta \in S^3 (w) \cap S^3(\zeta) $. 

Therefore, we have
\[
     f \circ \varphi_w (\theta)=g(\theta) {\mbox { and }} f(\theta)=g(\theta)  \;\;\; \forall \theta \in S^3(w) \cap S^3(\zeta) ,
\]
or
\[
     f \circ \varphi_w (\theta)=g(\theta) {\mbox { and }} f(\theta)=g(-\theta)  \;\;\; \forall \theta \in S^3(w) \cap S^3(\zeta).
\]
This implies that  either
\[
   f \circ \varphi_w (\theta)=f(\theta)  \;\;\; \forall \theta \in S^3(w) \cap S^3(\zeta)
\]
or
\[
   f \circ \varphi_w (\theta)=f(-\theta)  \;\;\; \forall \theta \in S^3(w) \cap S^3(\zeta),  
\]
where the restriction of $\varphi_w$ to $S^3(w) \cap S^3(\zeta)$ is not the identity, since we are assuming that $w \in S^3(\zeta)=\Xi_\alpha$. Thus, the restriction of $f$ to the 3 dimensional subspace spanned by $S^3(w) \cap S^3(\zeta)$ has an $O(3)$ symmetry, contradicting the hypothesis of Proposition \ref{funeq}. 

Since the case $S^3(\zeta)=\Xi_\alpha$ leads to a contradiction, we conclude  that either $S^3(\zeta)=\Xi_+$, or $S^3(\zeta)=\Xi_-$. Proposition \ref{funeq} is proven.

\qed

\section{Proof of Theorems \ref{tpr} and \ref{tse}} \label{sec4}

As in \cite{ACR}, the key ingredient in the proof of Theorem \ref{tpr} is the existence of a diameter $d_K(\zeta)$ such that the side projections of $K$ and $L$ are directly congruent.  We will first show that this implies that $L$ must also have a diameter in the $\zeta$ direction, which necessarily has the same length as  $d_K(\zeta)$. We can thus translate the bodies $K$ and $L$ so that their diameters $d_K(\zeta)$ and $d_L(\zeta)$ coincide and are centered around the origin. Since the translated bodies, $\tilde{K}$ and $\tilde{L}$, have countably many diameters, almost all side projections contain only this particular diameter, which must be fixed by the rotation. Therefore, we have reduced Theorem  \ref{tpr} to Proposition \ref{funeq} with $f=h_{\tilde{K}}$ and $g=h_{\tilde{L}}$.

\subsection{Theorem \ref{tpr} and Corollary \ref{cpr}.} 

Let $\zeta \in S^4$ be the direction of the diameter $d_K(\zeta)$ given in Theorem \ref{tpr}. By hypothesis, for every $w \in S^3(\zeta)$, the projections $K|w^\perp$ and $L|w^\perp$ are directly congruent. Hence, for every $w\in S^3(\zeta)$ there exists $\chi_w\in SO(4,S^3(w))$ and $a_w\in w^{\perp}$ such that
\begin{equation}\label{nax21}
\chi_w(K|w^{\perp})=L|w^{\perp}+a_w.
\end{equation}

%We will use  the following well-known properties of the support function. For every convex body $\tilde{K}$,  
%\begin{equation}\label{nax35}
%h_{\tilde{K}|w^{\perp}}(x)=h_{\tilde{K}}(x)   \;  \mbox{ and } \;  h_{\chi_w(\tilde{K}|w^{\perp})}(x)=h_{\tilde{K}|w^{\perp}}(\chi_w^t(x)), \quad \forall x \in w^{\perp},
%\end{equation}
%(see, for example, \cite[ (0.21), (0.26), pages 17--18]{Ga}).

Let ${\mathcal A}_{K}\subset S^4$ be the set of directions  parallel to the diameters of $K$, and ${\mathcal A}_{L}\subset S^4$ be the set of directions  parallel to the diameters of $L$. We define
\begin{equation}\label{nax3}
\Omega=\{w\in S^3(\zeta):\quad ({\mathcal A}_K\cup {\mathcal A}_L)\cap S^3(w)=\{\pm\zeta\}\}.
\end{equation} 

The following two Lemmata are proven by the same arguments used in \cite{ACR}. Lemma \ref{nunu2}  shows that for  most of the directions $w\in S^3(\zeta)$ the projections $K|w^{\perp}$ and $L|w^{\perp}$ have exactly one diameter, $d_K(\zeta)$ and $d_L(\zeta)$, respectively. We can thus translate the bodies $K$ and $L$ by vectors $a_K$, $a_L\in{\mathbb R}^5$,  to obtain $\tilde{K}=K+a_K$ and $\tilde{L}=L+a_L$ such that their diameters $d_{\tilde{K}}(\zeta)$ and $d_{\tilde{L}}(\zeta)$ coincide and are centered at the origin.

\bl(cf. \cite[Lemma 13]{ACR}.) \label{nunu2}
Let $K$ and $L$ be as in Theorem  \ref{tpr}, and let
$\zeta \in {\mathcal A}_K$. 
 Then    $\zeta \in {\mathcal A}_L$, and $\Omega$ is everywhere dense in $S^3(\zeta)$. Moreover, for every $w\in \Omega$ we have $\chi_w(\zeta)=\pm\zeta$ and $\omega_K(\zeta)=\omega_L(\zeta)$.
\el

\bl (cf. \cite[Lemma 14]{ACR}.) \label{naxuek1}
Let $\chi_w$ be the rotation given by (\ref{nax21}), and let $w\in \Omega$. Then the rotation $\varphi_w:=\left(\chi_w\right)^{t}$ %\in SO(3, S^2(w))$
satisfies $\varphi_w(\zeta)=\pm\zeta$ and 
\begin{equation}\label{naxuek2}
h_{\tilde{K}}\circ \varphi_w(\theta)=h_{\tilde{L}}(\theta)\qquad \forall\theta\in S^3(w).
\end{equation}
\el

\bigskip

\noindent{\bf Proof of Theorem \ref{tpr}.} Consider the closed sets  $\Xi=\{w\in S^3(\zeta):  (\ref{naxuek2}) \textrm{ holds with }\varphi_w(\zeta)=\zeta\}$ and $\Psi=\{w\in S^3(\zeta): (\ref{naxuek2}) \textrm{ holds with } \varphi_w(\zeta)=-\zeta\}$. Since the set $\Omega\subset(\Xi\cup \Psi)$ is everywhere dense in $S^3(\zeta)$ by   Lemma \ref{nunu2}, we have that $\Xi\cup \Psi= S^3(\zeta)$.  We have thus reduced matters to Proposition \ref{funeq} with $f=h_{\tilde{K}}$ and $g=h_{\tilde{L}}$.  
Therefore, either $h_{\tilde{K}}=h_{\tilde{L}}$ on $S^4$  or $h_{\tilde{K}}(\theta)=h_{\tilde{L}}(-\theta)$ for every $\theta\in S^4$. This means that either $K+a_K=L+a_L$ or  $K+a_K=-L-a_L$. 
\qed

\bigskip

\noindent{\bf Proof of Corollary \ref{cpr}.} First, we translate $K$ and $L$ by vectors $a_K, a_L \in \mathbb{R}^n$, obtaining the bodies $\tilde{K}=K+a_K$ and $\tilde{L}=L+a_L$, so that the diameters $d_{\tilde{K}}(\zeta)$ and $d_{\tilde{L}}(\zeta)$ are centered at the origin. Next, we observe that for any five dimensional subspace $J$ of $\mathbb{R}^n$, containing $\zeta$, the bodies $\tilde{K}|J$ and $\tilde{L}|J$ verify the hypotheses of Theorem \ref{tpr}. Therefore,  $\tilde{K}|J=\pm \tilde{L}|J$.

Assume that there exist two five dimensional subspaces $J_1$ and $J_2$, such that $\tilde{K}|J_1=\tilde{L}|J_1$ and $\tilde{K}|J_2=- \tilde{L}|J_2$. If $J_1 \cap J_2$ has dimension four, then 
$$\tilde{L}|(J_1\cap J_2)=(\tilde{L}|J_1)|(J_1\cap J_2)=(\tilde{K}|J_1)|(J_1\cap J_2)=(\tilde{K}|J_2)|(J_1\cap J_2)$$
\begin{equation}
 \label{minid}
=(-\tilde{L}|J_2)|(J_1 \cap J_2)=-\tilde{L}|(J_1\cap J_2).
\end{equation}
Since $-I \in SO(4)$, equation \eqref{minid} implies that the projection $\tilde{L}|(J_1\cap J_2)$ has an $SO(4)$ symmetry, contradicting the assumptions of the corollary. The same argument shows that if $ J_1 \cap J_2$ is three dimensional, the projection $\tilde{L}|(J_1\cap J_2)$ has an $O(3)$ symmetry (since $-I \in O(3)$). Next, assume that $ J_1 \cap J_2$ is two dimensional, and let $\{\zeta,v_1,v_2,v_3,v_4\}$ and $\{\zeta, v_1, v_2', v_3',v_4'\}$ be orthonormal bases for $J_1$ and $J_2$, respectively. %, where $v_2',v_3',v_4'$ are not in $J_1$. 
  Consider the subspace $J_0$ spanned by $\{\zeta, v_1, v_2,v_3,v_2'\}$.  Since $J_0$ is five dimensional and contains $\zeta$, we know that either $\tilde{K}|J_0=\tilde{L}|J_0$ or $\tilde{K}|J_0=-\tilde{L}|J_0$.  In the first scenario we have, 
$$\tilde{L}|(J_0\cap J_2)=(\tilde{L}|J_0)|(J_0\cap J_2)=(\tilde{K}|J_0)|(J_0\cap J_2)=(\tilde{K}|J_2)|(J_0\cap J_2)$$
\begin{equation}
 \label{minid2}
=(-\tilde{L}|J_2)|(J_0 \cap J_2)=-\tilde{L}|(J_0\cap J_2).
\end{equation}
Given that the dimension of $J_0 \cap J_2$ is 3, equation (\ref{minid2}) shows that $\tilde{L}|(J_0 \cap J_2)$ has an $O(3)$ symmetry, which is a contradiction.  For the second scenario, we apply the same argument to $J_0$ and $J_1$, showing that $\tilde{L}|J_0 \cap J_1$ has an $SO(4)$ symmetry (since the dimension of $J_1 \cap J_0$ is 4).  The case where $J_1\cap J_2$ is one dimensional can be dealt with in a similar way. %Let $\{\zeta,v_1,v_2,v_3,v_4\}$ and $\{\zeta, v_1', v_2', v_3',v_4'\}$ be orthonormal bases for $J_1$ and $J_2$, respectively.  We now consider $J_0=\{ \zeta, v_1, v_2, v_1', v_2'\}$.  Then the same argument as in the case where the dimension of $J_1 \cap J_2$ is 2 can be used.} 
 We conclude that  either $\tilde{K}|J=\tilde{L}|J$ for every five dimensional subspace $J$ containing $\zeta$, or $\tilde{K}|J=-\tilde{L}|J$ for all such $J$. 
By Theorem 3.1.1 from \cite[page 99]{Ga}, it follows that 
$\tilde{K}=\tilde{L}$ or  $\tilde{K}=-\tilde{L}$. Thus, 
 $K=L+a_L-a_K$ or $K=-L-a_L-a_K$.
 \qed

\section{Proofs of Theorem \ref{tse} and Corollary \ref{cse}.} \label{sec5}

We are now considering star-shaped bodies with respect to the origin.
Let  $\zeta\in S^4$ be the direction given in Theorem \ref{tse}. The hypotheses imply that 
%By the conditions of  Theorem \ref{ts},  the sections $K\cap w^{\perp}$ and $L\cap w^{\perp}$ are directly congruent for every $w\in S^2(\zeta)$.
for every $w\in S^3(\zeta)$ there exists $\chi_w\in SO(4,S^3(w))$ and $a_w\in w^{\perp}$ such that
\begin{equation}\label{nax211}
\chi_w(K\cap w^{\perp})=(L\cap w^{\perp})+a_w.
\end{equation}

Let $l(\zeta)$ denote the one dimensional subspace containing $\zeta$. As in Section \ref{sec3}, we let ${\mathcal A}_{K}\subset S^4$  be  the set of directions that are parallel to the diameters of $K$ (similarly for $L$). Note that it is possible for star-shaped bodies to contain several parallel diameters. We consider the set $\Omega$, defined as in \eqref{nax3}, and the set $\Omega^r$, defined by 
\begin{equation}
   \label{omr}
   \Omega^r=\{w\in \Omega: K\cap w^\perp {\mbox{ and }} L \cap w^\perp {\mbox{ have only one diameter}} \}.
\end{equation} 
Then it follows from the hypothesis of Theorem \ref{tse} that if $w \in \Omega^r$, then $d_K(\zeta)$ must be the unique diameter of $K \cap w^\perp$. We will use the notation $v_K(\zeta)=\rho_K(\zeta)+\rho_K(-\zeta)$ for the length of the diameter $d_K(\zeta)$. As in the previous Section, it can be shown that for most directions $w \in S^3(\zeta)$, the sections $K\cap w^\perp$, $L \cap w^\perp$ contain exactly one diameter parallel to $\zeta$ and passing through the origin.

\bl\label{radnu} (cf. \cite[Lemma 15]{ACR}.)
Let  $K$ and $L$ be  as in Theorem \ref{tse}.  Then $L$ has a diameter $ d_L(\zeta)$  passing through the origin,  and $\Omega^r$ is everywhere dense in $S^3(\zeta)$. Moreover, for every $w\in \Omega^r$ we have $\chi_w(\zeta)=\pm\zeta$  and $v_K(\zeta)=v_L(\zeta)$.
\el

%Since Lemma \ref{radnu} ensures that $L$ contains a diameter $d_L(\zeta)$ passing through the origin, 
We now wish to argue as in the proof of Theorem \ref{tpr} and translate the body $L$ so that its diameter $d_L(\zeta)$, given by Lemma \ref{radnu} coincides with $d_K(\zeta)$. However, the translate of a star-shaped body with respect to the origin may no longer be a star-shaped body with respect to the origin. The next Lemma, which is similar to Lemma 16 in \cite{ACR}  shows that, under our hypotheses, the translated body is still star-shaped (this Lemma is not necessary if $K$ and $L$ are convex). We include the proof here, since the rotation $\chi_w$ is now in  $SO(4,S^3(w))$, and the argument is slightly different.

\bl\label{stars}
There exists a vector $a \in \mathbb{R}^5$, parallel to $\zeta$, such that the body  $\tilde{L}=L+a$ is star-shaped with respect to the origin, and $d_K(\zeta)=d_{\tilde{L}}(\zeta)$. 
\el
\bp
Consider the sets
\begin{eqnarray*}
    R_1=\{ w \in \Omega^r: \; \chi_w(d_K(\zeta))=d_K(\zeta)\}, \\
    R_2=\{ w \in \Omega^r :\ \;  \chi_w(d_K(\zeta))\neq d_K(\zeta)\},
\end{eqnarray*}
where $\chi_w$ is the rotation in $SO(4,S^3(w))$ as in \eqref{nax211}. By Lemma \ref{rot}, since $\chi_w(\zeta)=\pm \zeta$, we have that either $\chi_w= \pm I$, or that $\chi_w$ has an invariant two dimensional subspace $\Pi^\perp$, containing $\zeta$, such that $\chi_w|\Pi^\perp=\pm I$. If $\chi_w=I$ or $\chi_w|\Pi^\perp= I$, or if $d_K(\zeta)$ is centered at the origin, then $w \in R_1$.  The only case in which $w \in R_2$ is if $d_K(\zeta)$ is not centered at the origin, and either $\chi_w=-I$ or $\chi_w|\Pi^\perp= -I$.

%If $w \in R_1$,  the rotation  must be either $\pm I$, or its a rotation about $\zeta$, or a rotation by angle $\pi$ about some $u \in S^2(\zeta)\cap S^2(w)$ (in this last case,  $d_K(\zeta)$ must be centered at the origin). On the other hand, if $w \in R_2$, $\chi_w$ is a rotation by angle $\pi$ about some $u \in S^2(\zeta)\cap S^2(w)$, and $d_K(\zeta)$ cannot be centered at the origin.

Assume, at first, that $\Omega^r = R_1$. Since the diameter $d_K(\zeta)$ is fixed by  $\chi_w$, and $d_L(\zeta)$ contains the origin, it follows that the vector $a_w$ in \eqref{nax211} is independent of $w\in \Omega^r$ and  $a_w=a_1=\left(\rho_K(\zeta)-\rho_L(\zeta) \right)\zeta$.
The translated section $(L\cap w^\perp )+a_1$ coincides with $\chi_w(K\cap w^\perp)$, and therefore $(L\cap w^\perp) +a_1$ is star-shaped with respect to the origin for every $w \in \Omega^r$. Since  $\Omega^r$ is dense in $S^3(\zeta)$, we conclude that the translated body $\tilde{L}=L+a$, with $a=a_1$, is also star-shaped with respect to the origin. 

Secondly, assume that  $\Omega^r = R_2$. Then,  $a_w$ is independent of $w \in \Omega^r$ and $a_w=a_2=
\left( \rho_K(-\zeta)-\rho_L(\zeta) \right) \zeta$. We conclude that $\tilde{L}=L+a$, with $a=a_2$, is star-shaped with respect to the origin.

Finally, we show that the case where $R_1$ and $R_2$ are both nonempty does not occur under the assumptions of Theorem \ref{tse}. %Let  $w_1 \in R_1$ and $w_2 \in  R_2$. 
Since $R_1 \cup R_2 =\Omega^r$, we have $S^3(\zeta)=\overline{R_1 \cup R_2} \subseteq \overline{R_1}\cup \overline{R_2} \subseteq S^3(\zeta)$.  Hence, there exists $w_0 \in \overline{R_1}\cap \overline{R_2}$, {\it i.e.}, there is a rotation $\chi_{w_0}$ such that $\chi_{w_0}(d_K(\zeta))=d_K(\zeta)$ and
\begin{equation}\label{r1}
   \chi_{w_0}(K\cap w_0^\perp)=L\cap w_0^\perp +a_1,
\end{equation}
and a rotation $\tilde{\chi}_{w_0}$ such that  $\tilde{\chi}_{w_0}(d_K(\zeta))\neq d_K(\zeta)$ and 
\begin{equation}\label{r2}
   \tilde{\chi}_{w_0}(K\cap w_0^\perp)=L\cap w_0^\perp +a_2.
\end{equation}
 In particular, since $\tilde{\chi}_{w_0}$ does not fix $d_K(\zeta)$, this diameter cannot be centered at the origin, and 
 it follows that  the other rotation $\chi_{w_0}$ must be the identity, at least on a two dimensional subspace containing $\zeta$. By \eqref{r1} and \eqref{r2} we have
\[
    K\cap w_0^\perp=\chi_{w_0}^{-1} \left( \tilde{\chi}_{w_0}(K\cap w_0^\perp)  \right) + b,
\] 
where $b \in \mathbb{R}^5$.
Observe that the rotation $\chi_{w_0}^{-1} \circ \tilde{\chi}_{w_0}$  is not the identity, since $\chi_{w_0}^{-1} \circ \tilde{\chi}_{w_0}(\zeta)=-\zeta$. Therefore, $K\cap w_0^\perp$ has a rotational symmetry. This contradicts the hypothesis of Theorem \ref{tse}. The Lemma is proven. 
\ep

In order to finish the argument, we need one further Lemma.

\bl\label{naxuxu1}(cf. \cite[Lemma 17]{ACR}.)
For every $w\in \Omega^r$ there exists $\varphi_w=\chi^{-1}_w\in SO(4, S^3(w))$, $\varphi_w(\zeta)=\pm\zeta$, such that
\begin{equation}\label{naxuxu2}
\rho_{K}\circ \varphi_w(\theta)=\rho_{\tilde{L}}(\theta)\qquad \forall \theta \in S^3(w).
\end{equation}
\el

\bigskip

\noindent{\bf Proof of Theorem \ref{tse}.} Consider the sets
$$
\Xi^r=\{w\in S^3(\zeta):\quad  (\ref{naxuxu2}) \quad\textrm{holds with}\quad \varphi_w(\zeta)=\zeta\}
$$
and 
$$
\Psi^r=\{w\in S^3(\zeta):\quad (\ref{naxuxu2}) \quad\textrm{holds with}\quad \varphi_w(\zeta)=-\zeta\}.
$$

By definition, $\Omega^r \subset (\Xi^r\cup \Psi^r)$. Therefore, Lemma  \ref{radnu} implies that  $\Xi^r\cup \Psi^r=S^3(\zeta)$. 
Now we can apply Proposition \ref{funeq} (with $f=\rho_{K}$, $g=\rho_{\tilde{L}}$, and $\Xi=\Xi^r$, $\Psi=\Psi^r$) 
obtaining that either $\rho_{K}=\rho_{\tilde{L}}$ on  $S^4$, or $\rho_{K}(\theta)=\rho_{\tilde{L}}(-\theta)$  for all $\theta\in S^4$. 
In the first case, $K=\tilde{L}$, and in the second, $K=-\tilde{L}$. Thus, either $K=L+a$, or $K=-L-a$.
This finishes the proof of Theorem \ref{tse}.
\qed

\bigskip

\noindent{\bf Proof of Corollary \ref{cse}}

The proof is similar to the one of Corollary \ref{cpr}. One has only to consider the sections $K \cap J$, $\tilde{L} \cap J$, instead of the projections $K|J$, $\tilde{L}|J$, and Theorem 7.1.1 from  \cite[page 270]{Ga}, instead of Theorem 3.1.1 from  \cite[page 99]{Ga}. 
\qed

\end{document}